\documentclass[10pt]{article}
\usepackage{amssymb}
\usepackage{latexsym,bm}
\usepackage{graphicx}
\usepackage{amsmath}
\usepackage{mathrsfs}

\date{}
\begin{document}
\title{A simple proof of Ore's theorem on the maximum size of $k$-connected graphs with given order and diameter\footnote{E-mail addresses:
{\tt mathdzhang@163.com}(L. Zhang).}}
\author{\hskip -10mm Leilei Zhang\thanks{Corresponding author.}\\
{\hskip -10mm \small Department of Mathematics, East China Normal University, Shanghai 200241, China}}\maketitle
\begin{abstract}
In 1968, Ore determined the maximum size of $k$-connected graphs with given order and diameter. We give a new short proof.
\end{abstract}


\vskip 8mm
We consider finite simple graphs and follow the book [3] for terminology and notations. Denote by $V(G)$ and $E(G)$ the vertex set and edge set of a graph $G,$ respectively. The {\it order} of a graph is its number of vertices, and the {\it size} is its number of edges. For a subset of vertices $S\subseteq V(G)$, we denote by $G[S]$ the subgraph of $G$ induced by $S$. The {\it diameter} of a graph $G$ is the greatest distance between two vertices of $G$. For a vertex $x$ in $V(G)$, the set $N_k(x)=\{y\in V(G):d_G(x,y)=k\}$ denotes the $k$-th neighborhood of $x$. Thus $N_1(x)$ is just the usual neighborhood $N(x).$
We use $G_1+\cdots+G_k$ to denote the {\it disjoint union} of graphs $G_1,\dots, G_k.$ The {\it join} of graphs $G_1$ and $G_2,$ written $G_1\vee G_2,$ is the graph obtained from $G_1+G_2$ by adding all the edges between $V(G_1)$ and $V(G_2).$ The {\it sequential join} of graphs $G_1,\dots,G_k,$ written $G_1\vee G_2\vee\cdots\vee G_k,$ is the graph obtained from $G_1+\cdots+G_k$ by adding all the edges between $V(G_i)$ and $V(G_{i+1}),$ $i=1,\ldots,k-1.$

Ore [1] determined the maximum size of a $k$-connected graph with given order and diameter, and characterized the corresponding extremal graphs. In 2019, Qiao and Zhan [2] gave a simple proof of Ore's theorem in the case $k=1.$ Using their ideas, we give a short simple proof of Ore's theorem for a general $k.$ Note that the problem is trivial when the diameter is equal to $1$ or $2.$

We denote by $|G|$ and $e(G)$ the order and size of a graph $G,$ respectively. For $A,\,B\subseteq V(G),$ the notation $[A,\, B]$ denotes the set of all edges with one endpoint in $A$ and the other endpoint in $B.$ For subgraphs $H$ and $R$ of $G,$ the notation $[H,\,R]$ means $[V(H),\,V(R)].$ We denote by $d(x,y)$
the distance between two vertices $x$ and $y.$

{\bf Theorem 1.} (Ore [1]) {\it Let $f(n,k,d)$ denote the maximum size of a simple $k$-connected graph of order $n$ and diameter $d$ with
 $d\ge 3.$ Then
$$
f(n,k,d)=\begin{cases}
\frac{(3d-5)k^2+(5-d)k}{2}+\binom{n-kd+k-2}{2}+3k(n-kd+k-2),\,\,{\rm if}\,\, d\ge4;\\
\frac{n^2-3n+2}{2},\,\,{\rm if}\,\,d=3.
\end{cases}
$$
This maximum size is attained by a graph $G$ if and only if $G$ is a sequential join of $d+1$ graphs
of the form $G_0\vee G_1\vee\cdots\vee G_d$ where each $G_i$ is a complete graph with $|G_0|=|G_d|=1$ and $|G_i|=k$
for all $1\le i\le d-1$ except for possibly one $i$ with $2\le i\le d-2$ such that $|G_i|>k$ or for possibly one $i$ with
$2\le i\le d-3$ such that $|G_i|>k$ and $|G_{i+1}|>k$ if $d\ge 4,$ and $|G_1|\ge k$ and $|G_2|\ge k$ if $d=3.$
}

{\bf Proof.} Let $G$ be a $k$-connected graph of order $n$ and diameter $d.$ Since $G$ has diameter $d$, there are vertices $x$ and $y$ which are at distance $d$. Denote $S_i=N_i(x),$ $i=0,1,\ldots,d.$ Then for $i=1,\ldots,d-1,$ each $S_i$ is a vertex cut of $G,$ which separates $x$ and $y.$
Since $G$ is $k$-connected, we have $|S_i|\ge k$ for $i=1,\ldots,d-1.$ Let $V_i\subseteq S_i$ with $|V_i|=k$ for
$i=1,\ldots,d-1$ and $V_0=\{x\},$ $V_d=\{y\}.$ Denote $H=G[\mathop{\cup}\limits_{i=0}^{d} V_i]$ and $R=G-V(H).$ Then $|H|=kd-k+2$ and $|R|=n-kd+k-2.$

Note that $e(G[V_i])\le \binom{k}{2}$ for $i=1,\ldots,d-1,$ $|[V_0,\,V_1]|=k,$
$|[V_j,\,V_{j+1}]|\le k^2$ for $1\le j\le d-2$ and $|[V_{d-1},\,V_d]|\le k.$ We have the following estimates
$$
e(H)\le \frac{(3d-5)k^2+(5-d)k}{2}, \eqno (1)
$$
$$
e(R)\le \binom{n-kd+k-2}{2}. \eqno (2)
$$
Let $z$ be a vertex of $R,$ let $p$ be the smallest subscript such that $z$ has a neighbor in $V_p$
and let $q$ be the largest subscript such that $z$ has a neighbor in $V_q.$ Let $u\in N(z)\cap V_p$ and let $v\in N(z)\cap V_q.$ Then $d(x,\,u)=p$ and $d(x,\,v)=q.$ If $q-p>2,$ then using the path $u,z,v$ we deduce that $q=d(x,v)\le d(x,u)+d(u,v)\le p+2<q,$ a contradiction. This shows that every vertex of $R$ can have neighbors in at most three consecutive $V_is.$ Consequently every vertex of $R$ can have at most $3k$ neighbors in $H.$ Suppose $d\ge 4.$ We have
$$
|[H,\,R]|\le 3k(n-kd+k-2).\eqno (3)
$$
Combining inequalities (1), (2) and (3) we obtain
$$
e(G)\le \frac{(3d-5)k^2+(5-d)k}{2}+\binom{n-kd+k-2}{2}+3k(n-kd+k-2). \eqno (4)
$$
Equality holds in (4) if and only if equality holds in each of (1), (2) and (3). If equality holds in (1), then every $V_i$ is a clique
and every vertex in $V_j$ is adjacent to every vertex in $V_{j+1}$ for $0\le j\le d-1;$
if equality holds in (2), then $R$ is a complete graph and the neighbors of vertices of $R$ in $H$ lie in at most four consecutive $V_is$ ; if equality holds in (3), then every vertex in $R$ has exactly $3k$ neighbors in three
consecutive $V_is.$ This shows that the maximum size is attained by $G$ if and only if $G$ is a graph of the form as stated in Theorem 1.

For the case $d=3,$ every vertex of $R$ can have at most $2k+1$ neighbors in $H.$ Thus
$$
|[H,\,R]|\le (2k+1)(n-2k-2).\eqno (5)
$$
Combining (1), (2) and (5) we obtain
$$
e(G)\le 2k^2+k+\binom{n-2k-2}{2}+(2k+1)(n-2k-2)=\frac{n^2-3n+2}{2}. \eqno (6)
$$
Equality holds in (6) if and only if equality holds in each of (1), (2) and (5), which force $G$ to be a graph of the stated form.
This completes the proof.   \hfill $\Box$

\vskip 5mm

{\bf Acknowledgement.} The author is grateful to Professor Xingzhi Zhan for his constant support and guidance. This research  was supported by the NSFC grant 12271170 and Science and Technology Commission of Shanghai Municipality (STCSM) grant 22DZ2229014.

{\bf Availability of data and materials.} Data sharing not applicable to this article as no datasets were generated or analysed during the current study.


\begin{thebibliography}{99}
\bibitem{1} O. Ore, Diameters in graphs. J. Combin. Theory 5(1968), 75-81.
\bibitem{2} P. Qiao and X. Zhan, The largest graphs with given order and diameter: a simple proof. Graphs Comb. 35(2019), 1715-1716.
\bibitem{3} D.B. West, Introduction to Graph Theory, Prentice Hall Inc., Upper Saddle River, 1996.
 \end{thebibliography}
\end{document}